\documentclass[12pt]{article}
\usepackage{amsmath,amssymb,latexsym,color}
\usepackage{graphicx}
\usepackage{pdfpages}
\begin{document}
\newcommand{\qbc}[2]{ {\left [{#1 \atop #2}\right ]}}
\newcommand{\anbc}[2]{{\left\langle {#1 \atop #2} \right\rangle}}
\newcommand{\be}{\begin{enumerate}}
\newcommand{\ee}{\end{enumerate}}
\newcommand{\beq}{\begin{equation}}
\newcommand{\eeq}{\end{equation}}
\newcommand{\bea}{\begin{eqnarray}}
\newcommand{\eea}{\end{eqnarray}}
\newcommand{\beas}{\begin{eqnarray*}}
\newcommand{\eeas}{\end{eqnarray*}}

\newcommand{\sn}{\mathfrak{S}_n}
\newcommand{\rsk}{\stackrel{\mathrm{RSK}}{\longrightarrow}}
\newcommand{\cc}{\mathbb{C}}
\newcommand{\rr}{\mathbb{R}}
\newcommand{\zz}{\mathbb{Z}}
\newcommand{\qq}{\mathbb{Q}}
\newcommand{\pp}{\mathbb{P}}
\newcommand{\nn}{\mathbb{N}}
\newcommand{\ff}{\mathbb{F}}
\newcommand{\ca}{\mathcal{A}}
\newcommand{\cp}{\mathcal{P}}
\newcommand{\mt}{M.I.T.}

\newcommand{\tr}{\textcolor{red}}
\newcommand{\tb}{\textcolor{blue}}
\newcommand{\tg}{\textcolor{green}}
\newcommand{\tm}{\textcolor{magenta}}
\newcommand{\tbn}{\textcolor{brown}}
\newcommand{\tp}{\textcolor{purple}}
\newcommand{\tn}{\textcolor{nice}}
\newcommand{\tor}{\textcolor{orange}}

\definecolor{brown}{cmyk}{0,0,.35,.65}
\definecolor{purple}{rgb}{.5,0,.5}
\definecolor{nice}{cmyk}{0,.5,.5,0}
\definecolor{orange}{cmyk}{0,.35,.65,0}

\begin{centering}
{\large\bf Enumerative and Algebraic Combinatorics in the 1960's and
  1970's}\\[1em]
{Richard P. Stanley}\\
{\small University of Miami}\\
{\small (version of 17 June 2021)}\\[3em]
\end{centering}

The period 1960--1979 was an exciting time for enumerative and
algebraic combinatorics (EAC). During this period EAC was transformed
into an independent subject which is even stronger and more active
today. I will not attempt a comprehensive analysis of the development
of EAC but rather focus on persons and topics that were relevant to my
own career. Thus the discussion will be partly autobiographical.

There were certainly deep and important results in EAC before
1960. Work related to tree enumeration (including the Matrix-Tree
theorem), partitions of integers (in particular, the Rogers-Ramanujan
identities), the Redfield-P\'olya theory of enumeration under group
action, and especially the representation theory of the symmetric
group, $\mathrm{GL}(n,\cc)$ and some related groups, featuring work by
Georg Frobenius (1849--1917), Alfred Young (1873--1940), and Issai
Schur (1875--1941), are some highlights. Much of this work was not
concerned with combinatorics \emph{per se}; rather, combinatorics was
the natural context for its development. For readers interested in the
development of EAC, as well as combinatorics in general, prior to
1960, see Biggs \cite{biggs}, Knuth \cite[{\S}7.2.1.7]{knuth3}, Stein
\cite{stein}, and Wilson and Watkins \cite{w-w}.

Before 1960 there are just a handful of mathematicians who did a
substantial amount of enumerative combinatorics. The most important
and influential of these is Percy Alexander MacMahon (1854-1929). He
was a highly original pioneer, whose work was not properly appreciated
during his lifetime except for his contributions to invariant theory
and integer partitions. Much of the work in EAC in the 60's and 70's
can trace its roots to MacMahon. Some salient examples are his
anticipation of the theory of $P$-partitions, his development of the
theory of plane partitions, and his work on permutation
enumeration. It is also interesting that he gave the Rogers-Ramanujan
identities a combinatorial interpretation.

William Tutte (1917--2002) was primarily a graph theorist, but in the
context of graph theory he made many important contributions to
EAC. He developed the theory of matroids, begun by Hassler Whitney
(1907--1989), into a serious independent subject.\footnote{One further
  early contributor to matroid theory is Richard Rado
  (1906--1989). For the early history of matroid theory, see
  Cunningham \cite{cunningham}.}  He was involved in the famous
project of squaring the square (partitioning a square into at least
two squares, all of different sizes) and in the enumeration of planar
graphs. He also played a significant role as a codebreaker during
World War II, and he remained mathematically active almost until his
death.

Three further mathematicians who were active in enumerative
combinatorics prior to 1960 (and afterwards) are John Riordan
(1903--1988), Leonard Carlitz (1907--1999), and Henry Gould
(1928--\ ). Although all three were quite prolific (especially
Carlitz), it would be fair to say that their work in combinatorics
cannot be compared to MacMahon's and Tutte's great
originality. Riordan wrote two books and over 70 papers on
EAC. Riordan also has the honor of originating the name ``Catalan
numbers'' \cite[pp.~186--187]{rs:cat}. He is to be admired for
achieving so much work in combinatorics, as well as some work in
queuing theory, without a Ph.D.\ degree. I met him once, at
Rockefeller University. Carlitz, whom I also met once (at Duke
University), was exceptionally prolific, with about 770 research
papers and 45 doctoral students in enumerative combinatorics and
number theory. His deepest work was in number theory but was not
appreciated until many years after publication. David Hayes
\cite{hayes} states that ``(t)his unfortunate circumstance is
sometimes attributed to the large number of his research papers.''
Gould has over 150 papers. His most interesting work from a historical
viewpoint is his collection \cite{gould1} of 500 binomial coefficient
summation identities. Nowadays almost all of them would be subsumed by
a handful of hypergeometric function identities, but Gould's
collection can nevertheless be useful for a non-specialist. Gould also
published bibliographies \cite{gould2} of Bell numbers and Catalan
numbers. A curious book \cite{gould3} of Gould originated from more
than 2100 handwritten pages of notes (edited by Jocelyn Quaintance) on
relating Stirling numbers of the first kind to Stirling numbers of the
second kind via Bernoulli numbers.

I will not attempt a further discussion of the many mathematicians
prior to 1960 whose work impinged on EAC. The ``modern'' era of EAC
began with the rediscovery of MacMahon, analogous to (though of course
at a much smaller scale) the rediscovery of ancient Greek mathematics
during the Renaissance.  Basil Gordon inaugurated (ignoring some minor
activity in the 1930's) a revival of the theory of plane
partitions. His first paper \cite{c-g} in this area (with
M.\,S.\ Cheema) was published in 1964, followed by a series of eight
papers, some in collaboration with his student Lorne Houten, during
the period 1968--1983. Gordon's contributions were quite substantial,
e.g., the generating function for symmetric plane partitions, but he
did miss what would now be considered the proper framework for dealing
with plane partitions, namely, $P$-partitions, the RSK algorithm, and
representation theory. These will be discussed later in this
paper. Gordon also began the ``modern era'' of the combinatorics of
integer partitions with his 1961 combinatorial generalization of the
Rogers-Ramanujan identities \cite{gordon}. Following soon in his
footsteps was George Andrews, whose 1964 Ph.D.\ thesis \cite{andrews}
marked the beginning of a long and influential career devoted
primarily to integer partitions, including the definitive text
\cite{andrews2}.

The 1960's and 1970's (and even earlier) also saw
some miscellaneous gems that applied linear algebra to combinatorics,
e.g., \cite{berlekamp}\cite{bose}\cite{g-p1}\cite{g-p2}\cite{lovasz},
paving the way to later more systematic and sophisticated
developments, and the development of spectral graph theory, e.g.,
Hoffman and Singleton \cite{h-s} (to pick just one random
highlight). For further applications of linear algebra to
combinatorics, see Matou\u{s}ek \cite{matousek}. For further early
work on spectral graph theory, see for instance the bibliographies of
Brouwer and Haemers \cite{b-h} and Cvetkovi\'c, Rowlinson, and Simi\'c
\cite{c-r-s}.

Another watershed moment in the modern history of EAC is the amazing
1961 paper \cite{schen} of Craige Schensted (1927--2021).
% April 12, 1927, in Mayfield ND (MO 382004)
% died January 22, 2021
Schensted was a mathematical physicist
with this sole paper on EAC. He was motivated by a preprint of a paper
on sorting theory (published as \cite{b-b}) by Robert Baer and Paul
Brock. Schensted defines the now-famous bijection between permutations
$w$ in the symmetric group $\sn$ and pairs $(P,Q)$ of standard Young
tableaux of the same shape $\lambda\vdash n$. (The notation
$\lambda\vdash n$ means that $\lambda$ is a partition of the
nonnegative integer $n$.) This bijection was extended to multiset
permutations by Knuth (discussed below) is now most commonly called
the RSK algorithm (or just RSK). The letter R in RSK refers to Gilbert
de Beauregard Robinson, who (with help from D.\,E.\ Littlewood) had
previously given a rather vague description of RSK. For further
details on the history of RSK see \cite[pp.~399--400]{ec2}. For
$w\in\sn$ we denote this bijection by $w\rsk (P,Q)$. Schensted proves
two fundamental properties of RSK: (a) the length of the longest
increasing subsequence of $w$ is equal to the length of the first row
of $P$ or $Q$, and (b) if $w=a_1 a_2\cdots a_n\rsk (P,Q)$, then the
``reverse'' permutation $a_n\cdots a_2 a_1$ is sent to
$(P^t,(Q^*)^t)$, where $^t$ denotes transpose and $Q^*$ is discussed
below. It is immediate from (a) and (b) that the length of the longest
decreasing subsequence of $w$ is equal to the length of the first
column of $P$.\footnote{Joel Spencer once informed me that he
  rediscovered the RSK algorithm and the results of Schensted in order
  to solve a problem \cite{rabinowitz} in the \emph{American
    Mathematical Monthly} proposed by Stanley Rabinowitz. Since
  Spencer's solution was not published he received no recognition,
  other than having his name listed as a solver, for this admirable
  accomplishment.}

It was Marcel-Paul Sch\"utzenberger (1920--1996) who first realized
that RSK was a remarkable algorithm that deserved further
study. Beginning in 1963 \cite{schutz} he developed many properties of
RSK and its applications to the representation theory of the symmetric
group, including the fundamental symmetry $w\rsk (P,Q)\Rightarrow
w^{-1}\rsk (Q,P)$, the combinatorial description of the tableau $Q^*$
defined above, and the theory of \emph{jeu de taquin}, a kind of
two-dimensional refinement of RSK. He realized that the formula
$Q^{**}=Q$ (obvious from the definition) can be extended to linear
extensions of any finite poset, leading to his theory of promotion and
evacuation \cite{schutz:evac}. I can remember being both mystified and
enthralled when I heard him lecture on this topic at M.I.T.\ on April
21, 1971. Eventually I became familiar enough with the subject to
write a survey \cite{rs:evac}. The writing style of Sch\"utzenberger
(and his later collaborator Alain Lascoux) could be very
opaque. Sch\"utzenberger once said (perhaps tongue in cheek) that he
deliberately wrote this way because mathematics should be learned
through struggle; it shouldn't be handed to someone on a silver
platter.\footnote{Recall Nietzsche's famous quote, ``Was ist
  Gl\"uck?---Das Gef\"uhl davon, da\ss\ die Macht w\"achst, da\ss\ ein
  Widerstand \"uberwunden wird.'' (Happiness is the feeling that power
  increases---that resistance is being overcome.)}  Beginning in 1978
\cite{l-s1} Sch\"utzenberger began a long and fruitful collaboration
with Alain Lascoux (1944--2013). Highlights of this work include the
plactic monoid and the theory of Schubert polynomials.

Dominique Foata (1934-- ) was perhaps the first modern researcher to look
seriously at the work of MacMahon on permutation enumeration. Foata's
first paper \cite{foata1} in this area appeared in 1963, followed in
1968 by his famous bijective proof \cite{foata2} of the
equidistribution of the number of inversions and the major index of a
permutation in $\sn$.  Foata and various collaborators did much
further work to advance enumerative combinatorics, involving such
areas as multiset permutations, tree enumeration, Eulerian
polynomials, rook theory, etc. Other researchers in the flourishing
French school of EAC before 1980 include Dominique Dumont, Jean
Fran\c{c}on, Germain Kreweras, Yves Poupard, Sch\"utzenberger, and
Xavier G\'erard Viennot. In particular, Kreweras (1918--1998)
\cite{kreweras} and Poupard \cite{poupard} launched the
far-reaching subject of noncrossing partitions.

The theory of parking functions is another thriving area of
present-day EAC. The first paper on parking functions \emph{per se}
was published in 1966 by Alan Gustave Konheim and Benjamin Weiss
\cite{k-w}, though the basic result that there are $(n+1)^{n-1}$
parking functions of length $n$ is equivalent to a special case of a
result of Ronald Pyke \cite[Lemma~1]{pyke} in 1959. Some further work
was done in the 1970's by Foata, Fran\c{c}on, Riordan and others, but
the subject did not really take off until the 1990's, when connections
were found with diagonal harmonics, symmetric functions, Lagrange
inversion, hyperplane arrangements, polytopes, Tutte polynomials,
noncrossing partitions, etc. See Yan \cite{yan} for a survey of much
of this more recent work.

A major impetus to the development of EAC was Gian-Carlo Rota
(1932--1999). He began his career in analysis but after a while was
seduced by the siren call of combinatorics. His first paper \cite{f-r}
in EAC (with Roberto Frucht) appeared in 1963 and was devoted to the
computation of the M\"obius function of the lattice of partitions of a
set. Rota was very ambitious and always wanted to see the ``big
picture.'' He realized that the M\"obius function of a (locally
finite) partially ordered set, first appearing in the 1930's in the
work of Philip Hall and Louis Weisner, had tremendous potential for
unifying much of enumerative combinatorics and connecting it with
other areas of mathematics. This vision led to his seminal paper
\cite{rota:f1} on M\"obius functions, with the rather audacious title
``On the foundations of combinatorial theory. I. Theory of M\"obius
functions.'' This paper had a tremendous influence on the development
of EAC and placing posets in a central role. For some further
information on Rota's influence on EAC, see the introductory essays by
Bogart, Chen, Goldman, and Crapo in \cite{rotacomb}\footnote{For my
  own account on how I was influenced by Foundations I and decided to
  work with Rota, see \cite{rs:ubcp}}.

Finite posets, despite their simple definition, have a remarkably rich
theory. In addition to the vast number of interesting examples and
special classes of posets, there are a surprising number of deep
results and questions that are applicable to \emph{all} finite
posets. These include incidence algebras and M\"obius functions,
M\"obius algebras, the connection with finite distributive lattices
\cite[{\S}3.4]{ec1}, $P$-partitions \cite[{\S}3.15]{ec1}, poset
polytopes \cite{rs:poly}\cite[Exer.~4.58]{ec1}, the order complex and
order homology \cite[{\S}3]{b-g-s}, Greene's theory \cite{greene2} of
chains and antichains, evacuation and promotion \cite{rs:evac},
correlation inequalities (in particular, the XYZ conjecture
\cite{shepp}), the $\frac 13$-$\frac 23$ conjecture (surveyed by
Brightwell \cite{brightwell}), the Chung-Fishburn-Graham conjecture on
heights of elements in linear extensions \cite[{\S}3]{rs:af},
dimension theory \cite{trotter}, etc. I should also mention the text
\cite{birkhoff} on lattice theory by Garrett Birkhoff
(1911--1996)\footnote{Three persons significant for combinatorics died
  in 1996: Sch\"utzenberger, Birkhoff, and P\'al (Paul) Erd\H{o}s.},
first published in 1940 with three editions altogether. Most of the
book deals with infinite lattices and was not combinatorial, but the
early chapters have some interesting combinatorial material on finite
lattices and posets. (The term ``poset'' is due to Birkhoff.)

%I like to think that posets satisfy the ``Goldilocks axioms'' for
%binary relations, just as groups satisfy the Goldilocks axioms for
%binary operations. A weakening of the group axioms to semigroups, for
%instance, gives too much freedom, while a strengthening, e.g., to
%abelian groups, gives too much structure. Of course semigroups and
%abelian groups are still worthwhile topics, but they lack the richness
%of groups. The situation for binary relations is similar. A good
%indication of the richness of finite poset theory is the many
%interesting [results,properties??]  that hold for \emph{all} finite
%posets, and the many interesting examples of specific posets. In the
%former [situation, category, class??], we have incidence algebras and
%M\"obius functions, M\"obius algebras, the connection with finite
%distributive lattices, Greene's theory \cite{greene} of chains and
%antichains, correlation inequalities (in particular, the $XYZ$
%conjecture \cite{shepp}), the 1/3-2/3 conjecture (surveyed by
%Brightwell \cite{brightwell}), the Chung-Fishburn-Graham conjecture on
%heights of elements in linear extensions \cite{rs:af}, etc.

The 1960's saw the development of an EAC infrastructure, in
particular, conferences, journals, textbooks, and prizes. The first
conference on matroid theory took place in 1964 (see
\cite{cunningham}). On January 1, 1967, the Faculty of Mathematics was
founded at the University of Waterloo (in Waterloo, Ontario, Canada),
which included a Department of Combinatorics and Optimization. To this
day it remains the only academic department in the world devoted to
combinatorics (though the Center for Combinatorics at Nankai
University functions somewhat similarly). The University of Waterloo
became a center for enumerative combinatorics with the arrival of
David Jackson in 1972, followed by Ian Goulden, who received his
Ph.D.\ from Jackson in 1979 and has a long and fruitful collaboration
with him. Much of their early work appears in their book \cite{g-j},
which is jam-packed with engaging results.

The Department of Combinatorics and Optimization at the University of
Waterloo sponsored three early conferences in combinatorics. The Third
Waterloo Conference on Combinatorics, held in 1968, was the only one
to have a Proceedings \cite{waterloo} and was the first combinatorics
conference that I attended. At that time I was a graduate student at
Harvard University. Most of the talks were on graph theory and design
theory. Jay Goldman and Rota discussed the number of subspaces of a
vector space over $\ff_q$ \cite{g-r}. For me it was a fantastic
opportunity to meet such legendary (to me, at any rate) mathematicians
as Elwyn Berlekamp (later to become my official host when I was a
Miller Fellow at U.\,C.\ Berkeley, 1971--73), Branko Gr\"unbaum,
William Tutte, Alan Hoffman, Crispin Nash-Williams, Johan Seidel,
Nicolaas de Bruijn, Richard Rado, Nathan Mendelsohn, Ralph Stanton,
Richard Guy, Frank Harary, Roberto Frucht, Claude Berge, Horst Sachs,
Denis Higgs, \emph{et al.} One incident sticks in my mind that
illustrates the nature of algebraic combinatorics at that time. During
an unsolved problem session, an expert on graph automorphisms
presented the conjecture that a vertex-transitive graph $\Gamma$ with
a prime number $p$ of vertices has an automorphism which cyclically
permutes all the vertices. I instantaneously saw (though I did not
mention it until after the problem session) that since a transitive
permutation group acting on an $n$-element $S$ set has a subgroup of
index $n$ (the subgroup fixing some element of $S$), the order of
Aut$(\Gamma)$ is divisible by $p$. Hence by elementary group theory
(Cauchy's theorem), Aut$(G)$ contains an element of order $p$, which
can only be a $p$-cycle. Rota later told me that the word had spread
about my proof and that many participants were impressed by it.

Another early conference (which I did not attend) that involved some
EAC was the Symposium in Pure Mathematics of the American Mathematical
Society, held in 1968 at UCLA.\footnote{Since the late 1940's UCLA has
  had a lot of combinatorial activity, as summarized by Bruce
  Rothschild \cite{roth}.} Of the 24 papers appearing in the
conference proceedings \cite{ucla}, six or so can be said to concern
EAC. A further early conference was the two-week June, 1969, meeting
in Calgary entitled Combinatorial Structures and Their Applications
\cite{calgary}, organized by Richard Guy and Eric Milner. Participants
related to EAC included David Barnette, Henry Crapo, Jack Edmonds
(who introduced the concept of a polymatroid at the meeting), Curtis
Greene, Branko Gr\"unbaum, Peter McMullen, John Moon, Leo Moser,
Tutte, and Dominic Welsh.

The University of North Carolina had some strength in combinatorics
headed by the statistician Raj Chandra Bose (1901--1987). Thomas
Dowling was a student of Bose who received his Ph.D.\ in 1967 and
stayed at UNC for several years afterwards. Bose and Dowling were the
organizers of a 1967 combinatorics conference at UNC. Of the 33 papers
in the proceedings \cite{unc1} of the 1967 conference, perhaps three
or four of them could be said to belong to EAC (by Henry Mann,
Riordan, Lentin-Sch\"utzenberger, and Hoffman). There were also some
papers on coding theory which had a strong EAC flavor, though nowadays
coding theory is rather tangential to mainstream EAC.

Bose and his collaborator K.\,R.\ Nair \cite{b-n} founded the theory
of \emph{association schemes} in 1939, though the term ``association
scheme'' is due to Bose and T. Shimamoto \cite{b-s} in 1952. They
arose in the applications of design theory to statistics. They became
of algebraic interest (primarily linear algebra) with a 1959 paper by
Bose and Dale Mesner \cite{b-m} and are a kind of combinatorial
analogue of the character theory of finite groups. A major
contribution was made in the 1973 doctoral thesis by Philippe Delsarte
at the Universit\'e Catholique de Louvain, reprinted as a Philips
Research Report \cite{delsarte}. Today association schemes remain an
active subject, but there is little interaction between its
practitioners and ``mainstream'' EAC.

A second, more elaborate combinatorics conference was held at UNC in
1970 (which I had the pleasure of attending). Its participants
included Martin Aigner, George Andrews, Edward Bender, Raj Chandra
Bose, Va\u{s}ek Chv\'atal, Henry Crapo, Philippe Delsarte, Jack
Edmonds, Paul Erd\H{o}s, Ray Fulkerson, Jean-Marie Goethals, Jay
Goldman, Solomon Golomb, Ralph Gomory, Ron Graham, Branko Gr\"unbaum,
Richard Guy, Larry Harper, Daniel Kleitman, Jesse MacWilliams, John
Moon, Ronald Mullin, Crispin St.\ John Alvah Nash-Williams, Albert
Nijenhuis, George P\'olya, Richard Rado, Herbert Ryser, Seymour
Sherman, Neil Sloane, Gustave Solomon, Joel Spencer, Alan Tucker, Neil
White, and Richard Wilson. Of the 53 papers in the Proceedings
\cite{unc2}, around 18 were in EAC and another four on the
boundary. (Of course there is some subjectivity in evaluating which
papers belong to EAC.)  One can see quite a significant increase in
EAC activity from 1967 to 1970. In particular, the 1970 proceedings
had a number of papers on matroid theory, a rapidly growing area
within EAC. Rota was paying some visits to UNC at this time and
managed to make many combinatorial converts, including Robert Davis,
Ladnor Geissinger, William Graves, and Douglas Kelly. Rota's student
Tom Brylawski also arrived at UNC in 1970.

We can mention two more meetings of note, both occurring in 1971.  The
University of Waterloo had a Conference on M\"obius Algebras. The
\emph{M\"obius algebra} of a finite lattice (or more generally, finite
meet-semilattice) over a field $K$ is the semigroup algebra over $K$
of $L$ with respect to the operation $\wedge$ (meet). It was first
defined by Louis Solomon, who developed its basic properties and
extended the definition to arbitrary finite posets. Davis and
Geissinger made further contributions, and Curtis Greene
\cite{greene:ma} gave an elegant presentation of the theory with
additional development, including a more natural and transparent
treatment of the extension to finite posets. An exposition (for
meet-semilattices) is given in \cite[{\S}3.9]{ec1}. M\"obius algebras
allow a unified algebraic treatment of the M\"obius function of a
finite lattice. One can see how much EAC progressed since Rota's
\emph{Foundations I} paper \cite{rota:f1} in 1964 by the existence of
a conference on the rather specialized topic of M\"obius algebras
(though the actual meeting was not really this specialized). For the
Proceedings, see \cite{mobalg}. One amusing aspect of this conference
is that, due to a tight budget, the participants from the USA shared a
single hotel suite. The suite had a large living room and two
adjoining bedrooms. Those unfortunate enough not to get a bedroom
slept in sleeping bags in the living room. Naturally one of the
bedrooms went to Rota. The second person (nameless here, but not
myself) received this honor because he snored loudly.

\begin{figure}
\centering
\centerline{\includegraphics[width=15cm]{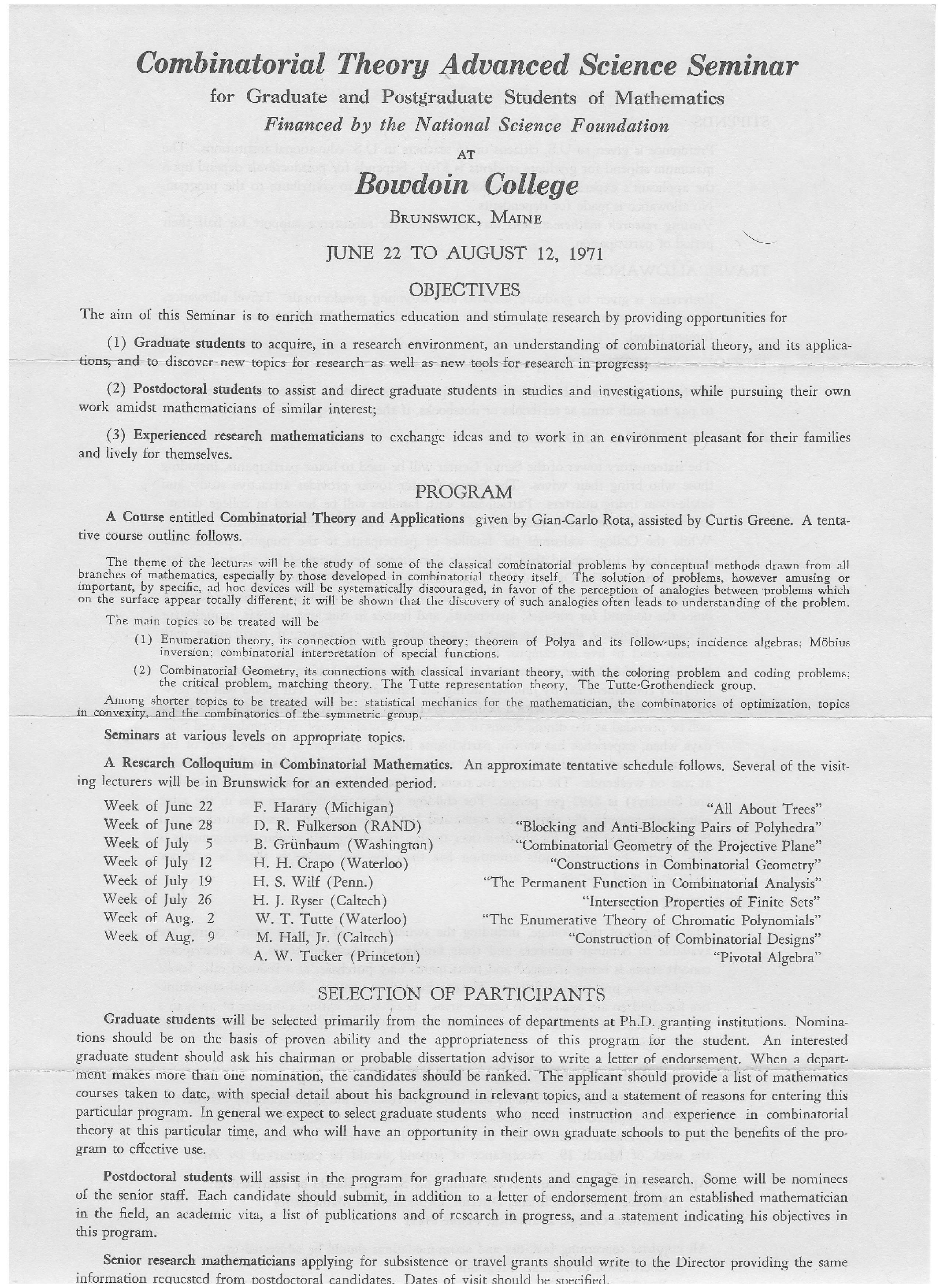}}
\caption{Bowdoin seminar announcement}
\label{fig2}
\end{figure}

The other 1971 meeting was the National Science Foundation sponsored
Advanced Seminar in Combinatorial Theory, held for eight weeks during
the summer at Bowdoin College in Brunswick, Maine. The existence of
this conference further exemplifies the tremendous progress made by
EAC since the early 1960's. Figure~\ref{fig2} shows an announcement of
the meeting. The main speaker was Rota, who gave a course on
combinatorial theory and its applications, assisted by Greene. There
were also nine featured speakers, one per week except two the last
week.

Rota and Greene soon decided that it would be better to alter their
lecture arrangement and for Greene to give instead a parallel set of
lectures on combinatorial geometries, with lecture notes written up by
Dan Kennedy. ``Combinatorial geometry'' was a term introduced by Rota
to replace ``matroid.''\footnote{More accurately, ``(combinatorial)
  pregeometry'' was the
    new term for matroid. A (combinatorial) geometry is a loopless
    matroid for which all two-element subsets are independent. Such
    matroids are called ``simple matroids.''}.  Rota thought that the
  term ``matroid'' was cacaphonous and had a frivolous
  connotation. However, Rota's terminology never caught on, so
  eventually he accepted the term ``matroid.'' But he was right on the
  button about the importance of matroids in EAC and in mathematics in
  general. His seminal book \cite{c-r} with Henry Crapo (Rota's first
  Ph.D.\ student in combinatorics) established the foundations of the
  subject for EAC.  Sadly only a preliminary edition was
  published. One of the topics of this book is the \emph{critical
    problem}, which Rota hoped would be a fruitful new approach toward
  the coloring of graphs. This hope, however, has not been realized.

We have mentioned Rota's great foresight in assessing the role of
posets and matroids within EAC. Two other subjects for which Rota had
high hopes are finite operator calculus (which began with a joint
paper \cite{m-r} with Ronald Mullin) and classical invariant
theory. Finite operator calculus can be appreciated by considering the
two operators $D$ and $\Delta$ on the space of all polynomials in one
variable $x$, say over the reals. They are defined by
 \beas Df(x) & = & f'(x)\ \ \mathrm{(differentiation)}\\
 \Delta f(x) & = & f(x+1)-f(x)\ \ \mathrm{(difference)}. \eeas
 They have many analogous properties; here are a sample. We use the
 falling factorial notation $(y)_m = y(y-1)\cdots (y-m+1)$.

 \begin{tabular}{c|c}
   $D$ & $\Delta$\\ \hline
   $De^x = e^x$ & $\Delta 2^x = 2^x$\\
   $Dx^n= nx^{n-1}$ & $\Delta (x)_n = n(x)_{n-1}$\\
   $f(x+t) = \sum_{n\geq 0} D^nf(t)\frac{x^n}{n!}$ &
   $f(x+t) = \sum_{n\geq 0} \Delta^n f(t)\frac{(x)_n}{n!}$.
 \end{tabular}

 Moreover, from the Taylor series expansion
  $$ f(x+1) = \sum_{n\geq 0}\frac{D^n f(x)}{n!} $$
 we get the formal identity $\Delta = e^D-1$. Finite operator calculus
 gives an elegant explanation for and vast generalization of these
 results. It also encompasses \emph{umbral calculus}, a seemingly
 magical technique made rigorous by Rota et al.\ in which powers like
 $a^n$ are replaced by $a_n$.\footnote{Gessel \cite{gessel} wrote a
 good paper for getting the flavor of umbral calculus.} The main
 shortcoming of finite operator calculus within EAC is its
 limited applicability.  Rota was a little miffed when I relegated
 finite operator calculus to a five-part exercise
 \cite[Exer.~5.37]{ec2} in my two books on enumerative combinatorics.

Adriano Garsia followed in Rota's footsteps by beginning in analysis
but converting to combinatorics (under Rota's influence). Garsia's
first combinatorics paper was an exposition \cite{garsia} of finite
operator calculus entitled ``An expos\'e of the Mullin-Rota theory of
polynomials of binomial type.'' Rota was annoyed at the title since
``expos\'e'' can have a negative connotation.\footnote{One dictionary
  definition is ``a report of facts about something, especially a
  journalistic report that reveals something scandalous.''} However,
Garsia simply intended for ``expos\'e'' to mean ``exposition.'' Garsia
went on to make many significant contributions to EAC, including (with
Stephen Milne) a long sought-for combinatorial proof of the
Rogers-Ramanujan identities \cite{g-m}. Later he did much beautiful
work related to symmetric functions, in particular, Macdonald
polynomials, diagonal harmonics, and the $n!$ conjecture.

Rota's work (with many collaborators) on classical invariant theory
never really caught on. One reason is that its algorithmic approach
toward algebraic questions (such as finite generation) were superseded
by the work of David Hilbert. Connections with algebraic geometry were
better handled by modern techniques, viz., the geometric invariant
theory pioneered by Mumford. More recently, however, there has been
somewhat of a resurgence of interest in classical invariant theory, as
attested by the book \cite{olver} of Peter Olver. Olver states that
this resurgence is ``driven by several factors: new theoretical
developments; a revival of computational methods coupled with powerful
new computer algebra packages; and a wealth of new applications,
ranging from number theory to geometry, physics to computer vision.''
Nevertheless, classical invariant theory has had limited impact on
modern EAC. There is scant mention of Rota in Olver's book.

My grades (using the common A--F grading system at most American
academic institutions) for the significance of Rota's four main
research topics within EAC are as follows:
 \begin{tabbing}aaaa\=aaaaaaaaaaaaaaaaaaaaaaaaaaaaaaa\=\kill
   \>{posets}\> {A+}\\ \>{matroid theory}\> {A}\\
   \>{finite operator calculus}\> {B--}\\
   \>{classical invariant theory}\> {C}
 \end{tabbing}
Matroid theory gets a slightly lower grade than posets because posets
are more ubiquitous objects both inside and outside EAC than
matroids. However, matroids have made many surprising appearances in
other areas that put them just below posets (in my opinion, of
course).\footnote{For information on the life and work of Rota, see
  \cite{kung}\cite{k-y}\cite{rota:mact}.}

One further topic of interest to Rota was Hopf algebras. In
particular, he anticipated the subject of combinatorial Hopf algebras,
e.g., in his perceptive article \cite{j-r} with Joni. I refrain from
giving a grade to Rota's work on Hopf algebras since he did not
develop the theory as intensively as the four topics above. For two
recent surveys of Hopf algebras in combinatorics, see Grinberg-Reiner
\cite{gr-re} and Loday-Ronco \cite{l-r}.

Getting back to the Bowdoin conference, let me mention three frivolous
anecdotes. One day Rota organized a tandem lecture. He chose six
people from the audience, including me. We had to leave the room and
not communicate among ourselves. He called us in the room one at a
time. Each person had to deliver a five-minute math lecture before the
next person was called in. The first lecture could be arbitrary, but
after that the lecture had to be a continuation of what was written on
the chalkboard. You were not allowed to erase anything that you wrote,
though you could erase what earlier speakers had written. My strategy
(not being the first person) was to decide on my lecture topic in
advance and to find some ridiculous way to link it to what appeared on
the board. It would be interesting to hold some more tandem lectures.

At the end of the meeting Rota also organized a prize ceremony. The
prizes were based on puns and other nonsense. I remember that I
received the Marriage Theorem (discussed by Rota in one of his
lectures) award because I arrived at the meeting two weeks late
shortly after my honeymoon. In addition to the prizes, Rota also
furnished several cases of champagne.

The final anecdote concerns what might be the all-time greatest EAC
pun. A famous theorem of Tutte characterizes graphic matroids in terms
of five excluded minors. Planar graphic matroids can be characterized
by two additional excluded minors (Kuratowski's theorem). Dan Kennedy
suggested in his writeup of Greene's talks on combinatorial geometries
that a \emph{pornographic matroid} is one for which all minors are
excluded.

A further indication of the growth of AEC was the occurrence of
workshops on this topic at the Mathematische Forschungsinstitut
Oberwolfach. It has held regular week-long workshops (eventually held
almost every week of the year) on different mathematical topics since
1949. The first workshop on combinatorics \emph{per se} was entitled
Angewandte Kombinatorik (Applied Combinatorics) and held October
15--21, 1972.\footnote{There were some previous workshops related to
  combinatorics, on such topics as discrete geometry, graph theory,
  and convex bodies.} All participants seem to have been
Europeans. Only a few were connected with EAC, including Foata,
Adalbert Kerber, and Volker Strehl.

\begin{figure}
\centering
\centerline{\includegraphics[width=17cm]{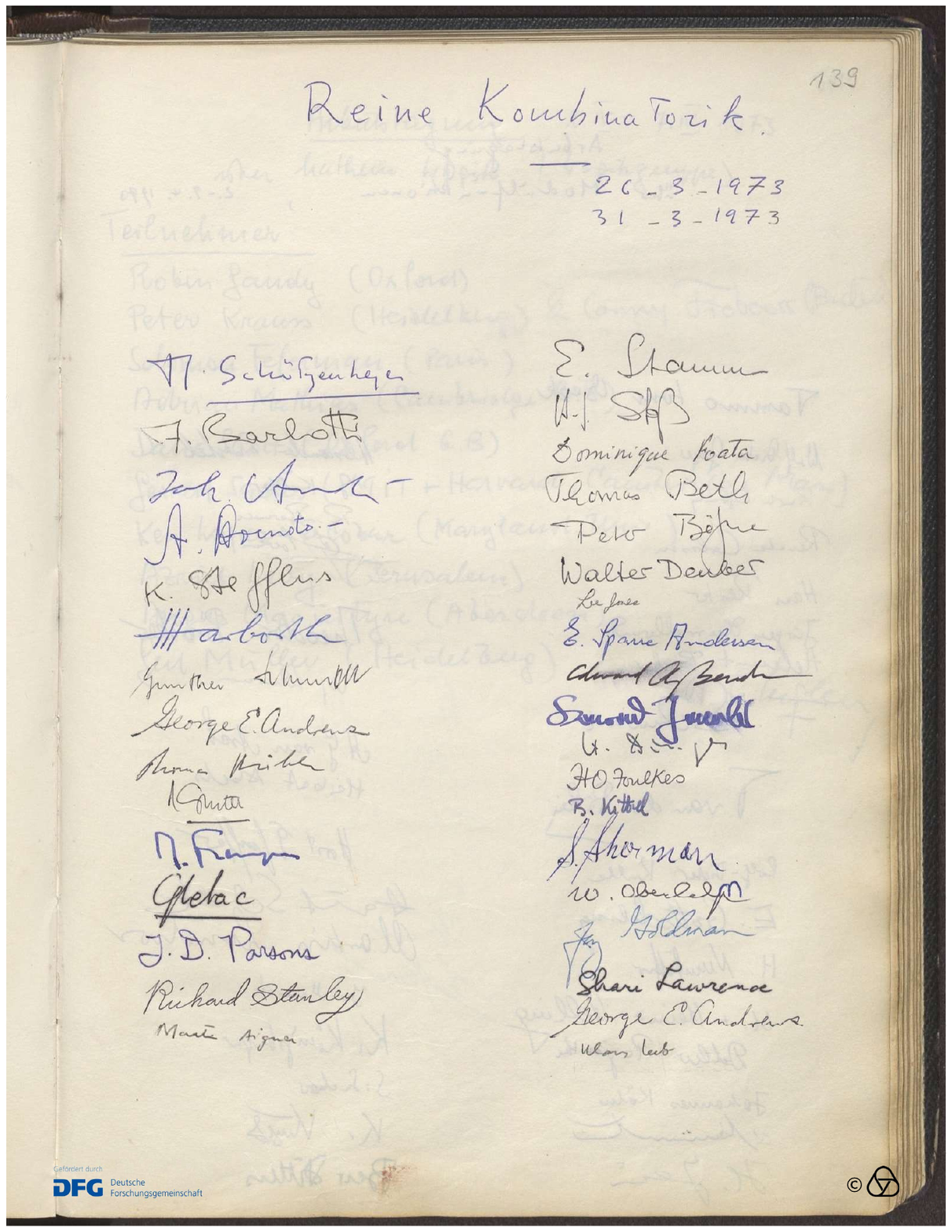}}
\caption{Oberwolfach participant list}
\label{fig:ober73}
\end{figure}

The second Oberwolfach workshop on combinatorics, entitled Reine
Kombinatorik (Pure Combinatorics), was held March 26--31, 1973 (my
first trip overseas). It was the first of several workshops organized
by Foata and Konrad Jacobs. The focus was almost entirely on EAC. See
Figure~\ref{fig:ober73} for the signatures of the participants
appearing in the \emph{G\"astebuch III}.\footnote{The Vortragsb\"ucher
  (Books of Abstracts) and G\"asteb\"ucher (Guest Books) of the
  Oberwolfach workshops are available online at the Oberwolfach Digital
  Archive \cite{oda}.} One highlight for me was the opportunity to
meet Herbert Foulkes, a pioneer of symmetric function combinatorics. I
gave a talk on my paper \cite{rs:ms}.

I had an amusing experience after this meeting unrelated to
mathematics. Martin Aigner invited me to visit him in T\"ubingen. I
then needed to take the train from T\"ubingen to the Frankfurt airport
for my return to the USA. It was necessary to change trains in
Stuttgart. Being unfamiliar with German cities and trains, when the
train arrived at Stuttgart-Bad Cannstatt I thought this was the main
station for Stuttgart and got off the train. By the time I realized my
error, the train had already departed. According to the posted
schedules, there was no way to get to the Frankfurt airport in time
for my flight. At the time I knew only a tiny bit of German and tried
to find someone in the station to whom to explain in English my
predicament. When that failed I went to a taxi stand and said that I
wanted to go to the Frankfurt airport. The driver was incredulous but
eventually said that the fare was around \$100 (US dollars), which was
a huge amount in 1973 for me to spend on a taxi. Nevertheless, I could
think of no feasible alternative so I agreed.  Given this sad
situation, who could believe that I would arrive at the Frankfurt
airport one hour before the train at a cost of \$5.00?

I showed the driver some German currency Travelers Checks (credit
cards were not in common use then) to show that I could pay. He was
not familiar with Travelers Checks and took me to an office in the
train station where they could be verified. Fortunately the woman in
the office could speak English. She told me that it might be cheaper
to fly from Stuttgart to Frankfurt. She called a travel agency and
found out that there was a special connector flight to my Frankfurt
flight, and that the cost was free! She got me a reservation, and the
taxi driver ended up driving me only to the Stuttgart airport for
around \$5. I ended up arriving at the Frankfurt airport about one
hour before the train! Had I known about this possibility in advance,
I could have even saved money by purchasing my train ticket only to
Stuttgart.

Conferences on EAC in the USA and Europe were becoming increasingly
more common. Another one in Oberwolfach (this time entitled just
Kombinatorik) in 1975 gave me the chance to meet for the first time
such luminaries as Kerber, Gilbert de B. Robinson, Gordon James,
Glanffrwd Thomas, Paul Stein, Louis Comtet, Ron King, Hanafi Farahat,
and Michael Peel. Many of these persons worked on symmetric functions
and the representation theory of the symmetric group, which was
quickly becoming a major area of EAC. I remember that Comtet gave
perfectly organized talks, with every square inch of the chalkboard
planned in advance, similar to later talks I heard by Ian
Macdonald.

One other conference in the mid-1970's was especially noteworthy for
me. This was the NATO Advanced Study Institute on Higher
Combinatorics, held in (West) Berlin on September 1--10, 1976
\cite{berlin}. A nonmathematical highlight of the meeting was a highly
structured and supervised visit to East Berlin. I remember that when
we returned to West Berlin, the East German guards used mirrors to
look under our bus to check for possible defectors. One of the
participants of the Berlin meeting was a graduate student from
Kungliga Tekniska h\"ogskolan (Royal Institute of Technology) in
Stockholm named Anders Bj\"orner. This meeting inspired him to work in
mainstream EAC, especially the emerging area of topological
combinatorics. He was a visiting graduate student at M.I.T.\ during
the academic year 1977--78 and went on to become a leading researcher
in EAC and a close collaborator of mine. In 1979 I was the opponent
(somewhat like an external examiner) for his thesis defense in
Stockholm.\footnote{A nice Swedish tradition is that the student
  treats the thesis defense attendees to dinner after the defense.} A
slightly earlier visitor to the Boston area was Louis Billera, who was
at Brandeis University for the 1974--75 academic year. He started out
in operations research and game theory but became interested in
combinatorial commutative algebra. At that time Brandeis was a world
center for commutative algebra, but Billera spent quite a bit of time
at M.I.T.\ and also developed into an EAC leader.

% Meeting related to EAC were becoming increasingly common. I will
% just mention one further one. banff?? 1981

%Another indication of the rapid expansion of EAC (and combinatorics in
%general) is the emergence of combinatorics textbooks, journals, and
%prizes.
Books on EAC saw a rapid development beginning in the 1960's. Before
then, we have Whitworth's treatise \emph{Choice and Chance} on
elementary enumeration and probability theory, first published in
1867. The 1901 fifth edition \cite{whitworth} includes 1000
exercises. Two further premodern books are the pioneering treatise
\cite{netto} by Eugen Netto, and the fascinating opus \cite{macmahon}
by MacMahon. Riordan's 1958 book \cite{riordan} is a kind of bridge
between the premodern and modern eras. Netto dealt primarily, and
MacMahon and Riordan exclusively, with enumerative combinatorics.

An interesting book \cite{d-b} by Florence Nightingale David and
D.\,E.\ Barton on enumerative combinatorics aimed at statisticians was
published in 1962. They followed up in 1966, joined by Maurice George
Kendall, with tables of symmetric function data (but with nothing on
Schur functions) preceded by a lengthy introduction \cite{d-k-b}.
Herbert Ryser in 1963 wrote an engaging monograph \cite{ryser} on
combinatorics containing some chapters on enumeration.  In 1967
Marshall Hall (my undergraduate adviser) wrote a textbook \cite{mhall}
that gave a quite broad coverage of combinatorics, including such
topics as permutations, M\"obius functions of posets, generating
functions, partitions, Ramsey theory, extremal problems, the simplex
method, de Bruijn sequences, block designs, and difference sets. There
appeared one year later a book \cite{berge} by Claude
Berge\footnote{Berge was primarily a graph theorist, but he was
  amazingly prescient in his selection of topics for his book on
  combinatorics.} (with a very entertaining Foreword by
Rota\footnote{This Foreword includes the famous statement, referring
  to the 1968 French edition,``Soon after that reading, I would be one
  of many who unknotted themselves from the tentacles of the Continuum
  and joined the then Rebel Army of the Discrete.''}) on EAC with many
interesting topics, including the M\"obius function of a poset and,
for the first time in a book, a discussion of standard Young tableaux,
RSK, and counting chains in Young's lattice. Also in 1969 there
appeared a book \cite{liu} by Chung-Laung Liu based on a course taught
in the Electrical Engineering Department of M.I.T. As in Hall's book
there is a broad selection of topics including some enumerative
combinatorics. In 1968 Riordan followed up his book \cite{riordan} on
combinatorial analysis with a book \cite{riordan2} on combinatorial
identities.

In 1964 appeared a collection of articles by leading mathematicians,
applied mathematicians, and physicists entitled \emph{Applied
  Combinatorial Mathematics} \cite{beck}, edited by Edwin
F. Beckenbach. The book were divided into four parts: Computation and
Evaluation, Counting and Enumeration, Control and Extremization, and
Construction and Existence. The part on Counting and Enumeration had
the articles ``Generating Functions'' by Riordan, ``Lattice
Statistics'' by Elliot W. Montroll, ``P\'olya's Theory of Counting''
by de Bruijn, and ``Combinatorial Problems in Graphical Enumeration''
by Frank Harary.  The authors of articles in the other three parts
include Derrick H. Lehmer, Richard Bellman, Albert W. Tucker, Marshall
Hall, Jr., Jacob Wolfowitz, and George Gamow. There are four
appendices reprinted from a 1949 English translation \cite{weyl}
(somewhat revised) of a 1926 article by Hermann Weyl. The first
appendix is entitled ``Ars Combinatoria.''

A very interesting set of lecture notes was published in 1969 by the
mathematical physicist Jerome Percus \cite{percus} (upgraded to a more
polished version \cite{percus2} in 1971). The first part deals with
enumerative topics like generating functions, MacMahon's Master
Theorem, partitions, permutations, and Redfield-P\'olya theory. The
second part concerns combinatorial problems arising from statistical
mechanics, such as the dimer problem, square ice, and the Ising
model. This deep topic is one of the most significant applications of
EAC to other areas and remains of great interest today. Let me just
mention the book \cite{baxter} of Rodney Baxter as an example of this
development.

A significant publication was Comtet's 1970 \emph{Analyse
  combinatoire} \cite{comtet1}, with an expanded English edition
\cite{comtet2} published in 1974. Despite the small physical size
($4\frac 12''\times 6\frac{15}{16}''$) of the 1970
French edition, it was packed with an incredible amount of interesting
enumerative combinatorics. The expanded English edition contains even
more information. Finally, mention should be made of Earl Glen
Whitehead's 1972 lecture notes \cite{whitehead}. Though rather
routine, it seems to be the first book or monograph with the title
``Enumerative Combinatorics'' or something similar.

According to MathSciNet, at the time of this writing (2021) there are
around 35 journals devoted to combinatorics or to the connection
between combinatorics and some other area. The first such journal did
not appear until 1966 (the year I started graduate school). This was
the \emph{Journal of Combinatorial Theory}, later split into parts A
(mainly EAC) and B (mainly graph theory). Frank Harary and Rota
founded the journal, Tutte was the first Editor-in-Chief, and Ron
Mullin became the de facto Managing Editor. The first issue had a
foreword by P\'olya, who called it ``a stepping stone to further
progress.'' Hans Rademacher, Robert Dickson, and Robert Plotkin had
the honor of having the first paper in this journal, followed by
Tutte. Further details on the founding of JCT and the subsequent split
into A and B are recounted by Edwin Beschler \cite{beschler} and the
editorial article \cite{50yrs}.

The first prize to be established in combinatorics was the George
P\'olya Prize in Applied Combinatorics, awarded by the Society for
Industrial and Applied Mathematics (SIAM). The recipients prior to
1980 were Ronald Graham, Klaus Leeb, Bruce Rothschild, Alfred Hales,
and Robert Jewett in 1971, Richard Stanley, Endre Szemer\'edi, and
Richard Wilson in 1975, and L\'aszl\'o Lov\'asz in 1979. The 1975
prize was awarded in San Francisco, so as a bonus Richard Wilson and I
(Szemer\'edi was not present) were invited by P\'olya, then aged 87,
to visit his home in Palo Alto. We spent an unforgettable evening
being shown by P\'olya his scrapbooks and other memorabilia. A further
combinatorics prize established prior to 1980 is the Delbert Ray
Fulkerson Prize of the American Mathematical Society. It was first
awarded in 1979, to Richard Karp, Kenneth Appel, Wolfgang Haken, and
Paul Seymour. Appel and Haken received it for their computer assisted
proof of the four-color conjecture.

One of my main mathematical interests, from graduate school to the
present day, is the theory of symmetric functions. I will briefly
discuss its rise to prominence in the 1960's and 70's. I already
mentioned the remarkable paper of Schensted \cite{schen} from 1961,
which was one of the main sources of stimulation for further
development of combinatorics related to symmetric functions. An even
more important paper for this purpose was published in 1959 in an
obscure conference proceedings by Philip Hall (1904--1982)
\cite{phall}.\footnote{Marshall Hall and Philip Hall are not related.}
The subject might have developed sooner if this publication had been
more widely known.\footnote{Around 1969 I gave a talk to my fellow
  graduate students at Harvard on Hall's paper. This was perhaps the
  first presentation of this material since Hall himself in 1959.} It
develops the now-familiar linear algebraic approach toward symmetric
functions, presenting known results in terms of five bases for
symmetric functions (monomial, elementary, complete, power sums,
Schur), the transition matrices between them, the involution $\omega$,
the scalar product making the Schur functions an orthonormal basis,
etc.  The MacTutor History of Mathematics \cite{phall-mt} gives the
following quote by Philip Hall:

\begin{quote}
 $\dots$ whenever in mathematics you meet with partitions, you have only
  to turn over the stone or lift up the bark, and you will, almost
  infallibly, find symmetric functions underneath. More precisely, if
  we have a class of mathematical objects which in a natural and
  significant way can be placed in one-to-one correspondence with the
  partitions, we must expect the internal structure of these objects
  and their relations to one another to involve sooner or later
  $\dots$ the algebra of symmetric functions.
\end{quote}

I learned of the paper of Philip Hall from Robert McEliece
(1942--2019), who was my colleague at the Caltech Jet Propulsion
Laboratory when I spent the summers of 1965--1970 there. McEliece had
studied for a year with Hall and had a copy of his paper
\cite{phall}. McEliece knew that I was becoming interested in
tableaux-like objects and showed me Hall's paper. I was immediately
transfixed by this beautiful connection between algebra and
combinatorics, and I realized that it would have many applications to
plane partitions. This resulted in my survey paper \cite{rs:pp}
entitled ``Theory and application of plane partitions.''  As I explain
in \cite{rs:papers}, my original title for this paper was ``Symmetric
functions and plane partitions.'' Rota had a keen political sense and
told me to change the title since he was angling for me to eventually
receive tenure in the Applied Mathematics Group of the
M.I.T.\ Department of Mathematics. When the Department of Mathematics
split into two groups in 1964 (discussed by Harvey
Greenspan in \cite[pp.~309--314]{greenspan}), the Pure Group was
not interested in combinatorics (represented by Rota and Daniel
Kleitman). The Applied Group was more accommodating, so combinatorics
ended up there. Any tenured faculty could choose whichever group they
wanted. Since I started in Applied and the other senior
combinatorialists were there, I stayed there throughout my
career. However, it would have been more interesting for me to be
involved with hiring and promotion decisions for algebraists and
number theorists rather than numerical analysts, PDE specialists,
etc. I should mention that the Pure Group now has a lot more respect
for combinatorics than it did in 1964.

There is one further underrecognized researcher on symmetric function
theory worthy of mention here. This is Dudley Ernest
Littlewood\footnote{Not to be confused with the better known John
  Edensor Littlewood, who in fact was D.\,E.\ Littlewood's tutor at
  Trinity College, Cambridge.} (1903--1979), who made many significant
contributions. These include some Schur function expansions of
infinite products, a product formula for the principal specialization
$s_\lambda(1,q,\dots, q^{n-1})$ which easily implies the
``hook-content formula'' \cite[Thm.~15.3]{rs:pp}, formulas for
super-Schur functions in a restricted number of variables, the
Jacobi-Trudi theorem for the orthogonal and symplectic groups, and the
expansion of the orthogonal and symplectic analogue of Schur functions
in terms of Schur functions. Much of his work is collected in his book
\cite{littlewood}. One reason his work was not adequately recognized
during his lifetime is his use of old-fashioned notation and
terminology.

Other mathematicians were becoming interested in symmetric functions
in the 1960's. Ronald Read wrote a mainly expository paper \cite{read}
on Redfield-P\'olya theory based on symmetric functions, in which
Schur functions play a prominent role.  We have already mentioned the
work of Sch\"utzenberger on RSK and his subsequent collaboration with
Lascoux. Sch\"utzenberger told Donald Knuth (1938-- ) about
Schensted's paper and its connection with the representation theory of
the symmetric group, thereby inspiring Knuth in 1970 to come up with a
generalization of RSK (and the notion of dual RSK, which for
permutations is the same as RSK) to permutations of multisets
\cite{knuth}. Lascoux pooh-poohed this work as straightforward and
unoriginal since it can be deduced easily from the original RSK (see
\cite[Lemma~7.11.6]{ec2}). In a strict sense Lascoux was correct, but
Knuth's version is necessary for a host of applications. Knuth himself
shows that it gives a combinatorial proof of the fundamental Cauchy
and dual Cauchy identities (though Knuth did not use this
terminology). Knuth also establishes the far-reaching ``Knuth
relations'' which determine when two permutations in $\sn$ have the
same insertion tableau under RSK. This result led to many further
developments, including Curtis Greene's important and subtle
characterization \cite{greene} of the shape of $P$ or $Q$ when $w\rsk
(P,Q)$ in terms of increasing and decreasing subsequences of $w$, and
the theory of the plactic monoid \cite{plactic} due to Lascoux and
Sch\"utzenberger. Knuth's combinatorial proofs of the Cauchy and dual
Cauchy identities were extended by William H. Burge \cite{burge}
to prove four similar identities due to D.\,E.\ Littlewood
\cite[second ed., p.~238]{littlewood}.

Two years after Knuth's paper, Edward Bender and Knuth showed that
Knuth's generalized RSK is the perfect tool for enumerating certain
classes of plane partitions, including plane partitions with at most
$r$ rows and with largest part at most $m$, and symmetric plane
partitions.\footnote{The argument for symmetric plane partitions
  easily extends to symmetric plane partitions with at most $r$ rows
  (and therefore contained in an $r\times r$ square), though Bender
  and Knuth do not mention this.}  Numerous other contributions to
enumerative combinatorics are scattered throughout Knuth's monumental
opus \emph{The Art of Computer Programming}, with the first volume
\cite{knuth2} appearing in 1968.

Another person bitten by the symmetric function bug was Ian
Macdonald. He was originally a leading algebraist. Some of his
algebraic results involved both combinatorics and representation
theory, such as his famous paper \cite{macd1} on affine root systems,
so he was in a good position to work on symmetric functions. His book
\cite{macd} was the first comprehensive treatment of the theory of
symmetric functions. In addition to the ``basic'' theory of
\cite[Ch.~7]{ec2}, it included such topics as polynomial
representations of GL$(n,\cc)$ (only outlined in \cite{ec2} without
proofs) in the setting of polynomial functors, characters of the
wreath product $G\wr \sn$ (where $G$ is any finite group), Hall
polynomials, Hall-Littlewood symmetric functions, and the characters
of GL$(n,\ff_q)$. One highlight for combinatorialists are some
``examples'' \cite[Exam.~I.5.13--19]{macd,macd2} which use root
systems to unify many results and conjectures on plane partitions.  (A
vastly expanded second edition \cite{macd2}, with contributions from
Andrey Zelevinsky, contains, among other things, a discussion of what
are now called Macdonald polynomials.)  For a combinatorialist, the
biggest defect of this book (either edition) is the omission of
RSK. In fact, Macdonald once told me that his main regret regarding
his book was this omission.

A major area of EAC that developed in the 1960's and 70's was its
connections with geometry. The three main topics (all interrelated)
discussed here are (1) simplicial complexes, (2) hyperplane
arrangements, and (3) Ehrhart theory. Of course there was already lots
of work on connections between combinatorics and geometry prior to
1960. Much of this work (such as Helly's theorem) belongs more to
extremal combinatorics than EAC so will not be considered here. In
regard to simplicial complexes, the main question of interest here is
what can be said about the number of $i$-dimensional faces of a
simplicial complex\footnote{All simplicial complexes considered here
  are finite abstract simplicial complexes.} satisfying certain
properties, such as being pure (all maximal faces have the same
dimension), triangulating a sphere, having vanishing reduced homology,
etc. If a $(d-1)$-dimensional simplicial complex $\Delta$ has $f_i$
faces of dimension $i$ (or cardinality $i+1$), then the
\emph{$f$-vector} of $\Delta$ is defined by
$f(\Delta)=(f_0,f_1,\dots,f_{d-1})$.  The first significant
combinatorial result concerning $f$-vectors is the famous
Kruskal-Katona theorem, which characterizes $f$-vectors of arbitrary
simplicial complexes. This result was first stated without proof by
Sch\"utzenberger \cite{schutz:kk} in 1959, and then rediscovered
independently by Joseph Kruskal \cite{kruskal} and Gyula Katona
\cite{katona}.  It turns out that Francis Macaulay
\cite{mac2}\cite{mac} had previously given a multiset analogue of the
Kruskal-Katona theorem which he used to characterize Hilbert functions
of graded algebras. This result played an important role in later EAC
developments. A common generalization of the Kruskal-Katona theorem
and the Macaulay theorem is due to George Clements and Bernd
Lindstr\"om \cite{c-l}.  Further progress on $f$-vectors of simplicial
complexes was greatly facilitated by the introduction of tools from
commutative algebra, and later, algebraic geometry.

Around 1975 Melvin Hochster and I independently defined a ring (in
fact, a graded algebra over a field $K$) $K[\Delta]$ associated with a
simplicial complex $\Delta$. Hochster was interested in algebraic and
homological properties of this ring and gave it to his student Gerald
Reisner for a thesis topic. Hochster was motivated by his paper
\cite{hochster1}, a pioneering work in the interaction between
commutative algebra and convex polytopes, while I was interested in
the connection between $K[\Delta]$ and the $f$-vector of $\Delta$. The
paper \cite{reisner} of Reisner (based on his Ph.D.\ thesis) was just
what was needed to obtain information on $f$-vectors, beginning with
\cite{rs:con}\cite{rs:ubc}.
%For further historical information see \cite{rs:ubcp}.

As I explain in \cite{rs:ubcp}, I was led to define the ring
$K[\Delta]$ (which I denoted $R_\Delta$) by my previous work
\cite{rs:ms} on ``magic squares.'' Specifically, let $H_n(r)$ denote
the number of $n\times n$ matrices of nonnegative integers for which
every row and column sums to $r$. Anand, Dumir, and Gupta \cite{a-d-g}
conjectured that $H_n(r)$ is a polynomial in $r$ of degree $(n-1)^2$
with certain additional properties. I proved polynomiality by showing
that $H_n(r)$ is the Hilbert polynomial of a certain ring. This was
the beginning of the subject of \emph{combinatorial commutative
  algebra} \cite{miller}\cite{rs:ca}, now a well-established
constituent of EAC. Just after the paper \cite{rs:ms} was written, I
became aware that $H_n(r)$ is also the Ehrhart polynomial of the
Birkhoff polytope of $n\times n$ doubly stochastic matrices, leading
to a more elementary proof of the polynomiality of $H_n(r)$ (and some
other properties). Ehrhart theory is discussed below.

Algebraic geometry had a long history of connections with
combinatorics, such as the Schubert calculus and determinantal
varieties. In the 1970's the theory of toric varieties
\cite{danilov}\cite{k-k-m-s} established connections with convex
polytopes and related objects. Applications of algebraic geometry to
combinatorics first appeared in 1980 \cite{rs:hl}\cite{rs:g}, a bit
too late for this paper. Subsequently there developed the flourishing
subject \emph{combinatorial algebraic geometry}, with many related
papers, books, and conferences.

An important subject within EAC is the theory of hyperplane
arrangements, or just \emph{arrangements}, i.e., a discrete (usually
finite) collection of hyperplanes in a finite-dimensional vector space
over some field. If we remove a finite collection $\ca$ of hyperplanes
in $\rr^d$ from $\rr^d$, we obtain a finite number of open connected
sets called \emph{regions}. The closure of each region is a convex
polyhedron $\cp$ (the intersection, possibly unbounded, of finitely
many closed half-spaces). A \emph{face} of $\ca$ is a relatively open
face of one of the polyhedra $\cp$. In particular, a region is a
$d$-dimensional face. The primary combinatorial problem concerning
$\ca$ is the counting of its faces of each dimension $k$, and
especially the number of regions. There was some early work on special
classes of arrangements, such as those in $\rr^3$ and those whose
hyperplanes are in general position. Robert Winder \cite{winder} in
1966 was the first to give a general result---a formula for the number
of regions for any $\ca$ when all the hyperplanes contained the
origin.

Arrangements did not get established as a separate subject until the
pioneering 1974 thesis of Thomas Zaslavsky \cite{zas1}\cite{zas2}, one
of the most influential Ph.D.\ theses in EAC. He gave formulas for the
number of $k$-faces and bounded $k$-faces in terms of the M\"obius
function of the \emph{intersection poset} $L_\ca$ of all nonempty
intersections of the hyperplanes in $\ca$. Not only did this work
initiate the subject of hyperplane arrangements within EAC, but also
it cemented the role of the M\"obius function as a fundamental EAC
tool. If the hyperplanes of $\ca$ all contain the origin, then $L_\ca$
is a \emph{geometric lattice}. Geometric lattices are equivalent to
simple matroids, so Zaslavsky's work also elucidated the connection
between arrangements and matroids. This connection led to the theory
of oriented matroids, an abstraction of real arrangements (where one
can consider on which side of a hyperplane a face lies) analogous to
how matroids are an abstraction of linear independence over any
field. The original development (foreseen by Ralph Rockafellar
\cite{rock} in 1967) goes back to Robert Bland \cite{bland} in 1974,
Jim Lawrence \cite{lawrence} in 1975, and Michel Las Vergnas
\cite{lasv} in 1975.  Published details of this work appear in Robert
Bland and Las Vergnas \cite{b-l} and Jon Folkman and Jim Lawrence
\cite{f-l}, both in 1978.

The last area of EAC to be discussed here is the combinatorics of
integral (or more generally, rational) convex polytopes in $\rr^n$,
i.e., convex polytopes whose vertices lie in $\zz^n$ (or $\qq^n$). The
first inkling of the general theory is the famous formula of Georg
Pick (1859--1942) \cite{pick} in 1899 for the area $A$ enclosed by a
simple (i.e., not self-intersecting) polygon $\cp$ in $\rr^2$ with
integer vertices in terms of the number of the number $i$ of integer
points in the interior of $\cp$ and the number $b$ of integer points
on the boundary (that is, on the polygon itself), namely, $A=i+\frac
b2 -1$. The question naturally arises of extending this result to
higher dimensions.  John Reeve did some work for three dimensions in
the 1950's.  The first general result is due to Macdonald
\cite{macd:lp} in 1963.

In the meantime Eug\`ene Ehrhart was slowly developing since the
1950's his theory of integer points in polyhedra and their dilations,
including his famous ``loi de reciprocit\'e,'' in a long series of
papers published mostly in the journal \emph{Comptes Rendus
  Math\'ematique}. He gave an exposition of this work in a monograph
\cite{ehrhart} of 1974.\footnote{Ehrhart and his followers deal with
  \emph{convex} polytopes, while Pick's theorem holds for arbitrary
  simple polygons. The polytopal Ehrhart theory can be extended to
  integral (or rational) polyhedral complexes $\Gamma$, e.g.,
  \cite[{\S}4]{rs:lc}.} Ehrhart's trailblazing work was not
completely accurate. A rigorous treatment was first given by Macdonald
\cite{macd:cc} in 1971. Ehrhart theory is now a central area of EAC
with many applications and connections to other subjects.

\begin{figure}
\centering
\centerline{\includegraphics[width=6cm]{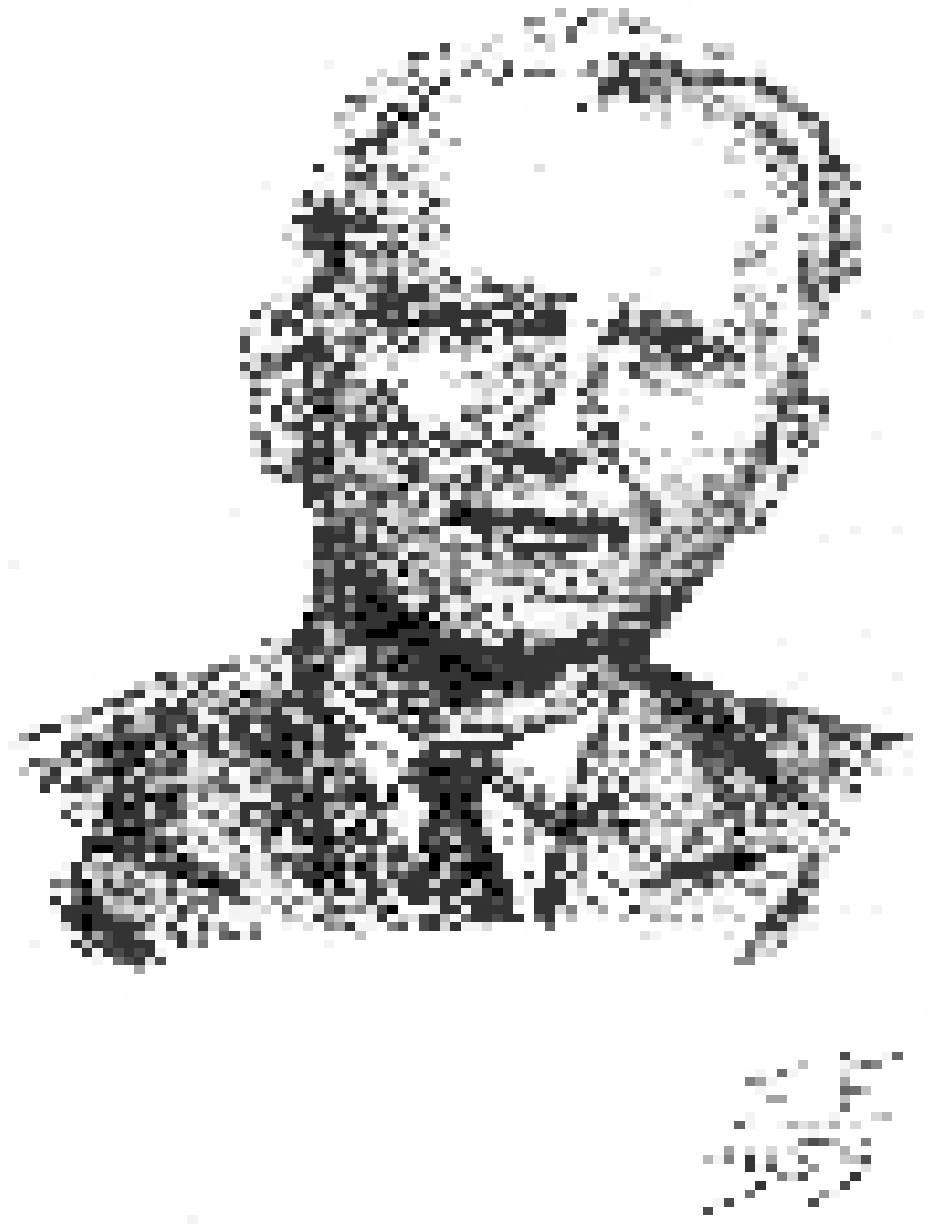}}
\caption{Eug\`ene Ehrhart}
\label{fig:ehrhart}
\end{figure}

 Eug\`ene Ehrhart (1906-2000) had an interesting career
 \cite{clauss}. In particular, he was a teacher in various lyc\'ees
 (French high schools) and did not receive his Ph.D.\ degree until the
 age of 59 or 60. He spent the latter part of his life in
 Strasbourg. Most of the professors in the Mathematics Department of
 the University of Strasbourg regarded him as somewhat of an amateur
 crackpot. When I gave a talk related to Ehrhart theory at this
 Department, the audience members were talking incredulously about
 whether the Ehrhart of Ehrhart theory was the same Ehrhart they
 knew. Ehrhart also had some artistic talent; Figure~\ref{fig:ehrhart}
 is a self-portrait.

 I had the pleasure of meeting Ehrhart once, at a 1976 conference
 ``Combinatoire et repr\'esentation du group sym\'etrique'' organized
 by Foata in Strasbourg \cite{crgs}. Some highlights of this
 conference are a survey by Robinson of the work of Alfred Young,
 Viennot's geometric version of RSK, and Sch\"utzenberger's
 presentation of the first proof of the validity of \emph{jeu de
   taquin}. Macdonald was in attendance but did not speak, certainly
 the last time such a state of affairs occurred at a combinatorics
 conference!

 At this conference Curtis Greene and I also had the pleasure of
 meeting Bernard Morin (1931--2018). Morin was a member of the first
 group to exhibit explicitly a crease-free eversion (turning
 inside-out) of the 2-sphere, and he also discovered the \emph{Morin
   surface}, a half-way point for the sphere
 eversion.\footnote{Stephen Smale was the first to prove that a
   crease-free eversion exists, but he did not give an explicit
   description. When Smale, as a graduate student at the University of
   Michigan, told his thesis adviser Raoul Bott that he (Smale) had
   proved that the 2-sphere can be everted, Bott replied that Smale
   must have made a mistake! Bott had in mind a simple argument
   showing the impossiblility of eversion, but Bott's reasoning was
   faulty.} Morin explained the sphere eversion in his office to
 Curtis Greene and me with the use of some models. The remarkable
 aspect of this story is that Morin was blind since the age of six!

 I will conclude this paper with a discussion of EAC at M.I.T.\ in the
 1960's and 1970's. Thanks to the influence of Rota, especially after
 he returned to \mt\ from Rockefeller University in 1967,
 \mt\ became a leading center for EAC research. His return to
 \mt\ conveniently coincided with the time when I was getting
 interested in some combinatorial problems as a Harvard graduate
 student, so I could experience the development of EAC at \mt\ almost
 from the beginning.

  It was certainly a thrilling time to be at \mt\ and interact with a
  plethora of graduate students, postdocs, visitors, seminar speakers,
  and the two senior combinatorialists, Kleitman and Rota.  Rota's
  first graduate student in combinatorics at \mt\ was Henry Crapo
  (1932--2019), who received his degree related to matroid theory in
  1964.\footnote{Among Rota's pre-combinatorics graduate students was
    Peter Duren. Peter Duren was a secondary thesis adviser of
    Theodore Kaczynski, the notorious Unibomber. Thus I am a
    ``secondary academic uncle'' of Kaczynski.} Rota's other
  combinatorics students in the 1960's and 1970's included Thomas
  Brylawski (from Dartmouth), Thorkell Helgason (later holding several
  prominent government positions in Iceland), Walter Whiteley, Stephen
  Fisk (Harvard), Neil White (Harvard), Peter Doubilet, Stephen Tanny,
  Kenneth Holladay, Hien Nguyen, Joseph Kung, and Joel Stein
  (Harvard). The combinatorics instructors during the late 1960's
  consisted of Curtis Greene, Michael Krieger, and Bruce Rothschild.

 Intersecting my time as a graduate student at Harvard were Edward
 Bender and Jay Goldman. Bender was a Benjamin Peirce Instructor, and
 Jay Goldman was a junior faculty member in the Department of
 Statistics. Bender worked primarily in what today is called
 \emph{analytic combinatorics}, though I mentioned above his role in
 the enumeration of plane partitions. Jay Goldman was trained as a
 statistician but was lured to combinatorics by Rota's spell. In the
 1970's he wrote five influential papers with James Joichi and Dennis
 White which established the subject of \emph{rook theory}.

 During the academic year 1967--68 I took from Bender and Goldman a
 graduate course in combinatorics, my
 first course on this subject.\footnote{Although Caltech was a center
   for combinatorics when I was an undergraduate there, at that time I
   did not think that combinatorics was a serious subject and declined
   to take any courses in this area!} Officially, the Bender-Goldman
 course was entitled ``Statistics 210: Combinatorial Analysis'' in the
 fall of 1967, taught by Goldman, and ``Mathematics 240: Combinatorial
 Analysis'' in the spring of 1968, taught by Bender. This course was
 the  second course in combinatorics offered at Harvard, the first
 being  ``Mathematics 240: Combinatorial Analysis'' taught by Alfred
 Hales from Ryser's book \cite{ryser} in the fall of
 1965.\footnote{There were earlier undergraduate courses on Applied
   Discrete Mathematics in the Applied Mathematics Department, but
   these were not really courses in combinatorics.} There was some
 interest around the time of the
 Bender-Goldman course in giving a unified development of generating
 functions. How to explain why generating functions like $\sum
 f(n)x^n$ and $\sum f(n)\frac{x^n}{n!}$ occurred frequently, while one
 never saw $\sum f(n)\frac{x^n}{n^2+1}$, for instance? Three theories
 soon emerged: (1) binomial posets, due to Rota and me
 \cite[{\S}8]{d-r-s}\cite[{\S}3.18]{ec1}, (2) dissects, due to Michael
 Henle \cite{henle}, and (3) prefabs, due to Bender and Goldman
 \cite{b-g}. The theory of prefabs was part of Bender and Goldman's
 course. None of these theories have played much of a role in
 subsequent EAC developments because of their limited
 applicability. Later, Andr\'e Joyal developed the theory of species
 \cite{joyal}, based on category theory, which is probably the
 definitive way to unify generating functions.
  
\begin{figure}
\centering
\centerline{\includegraphics[width=12cm]{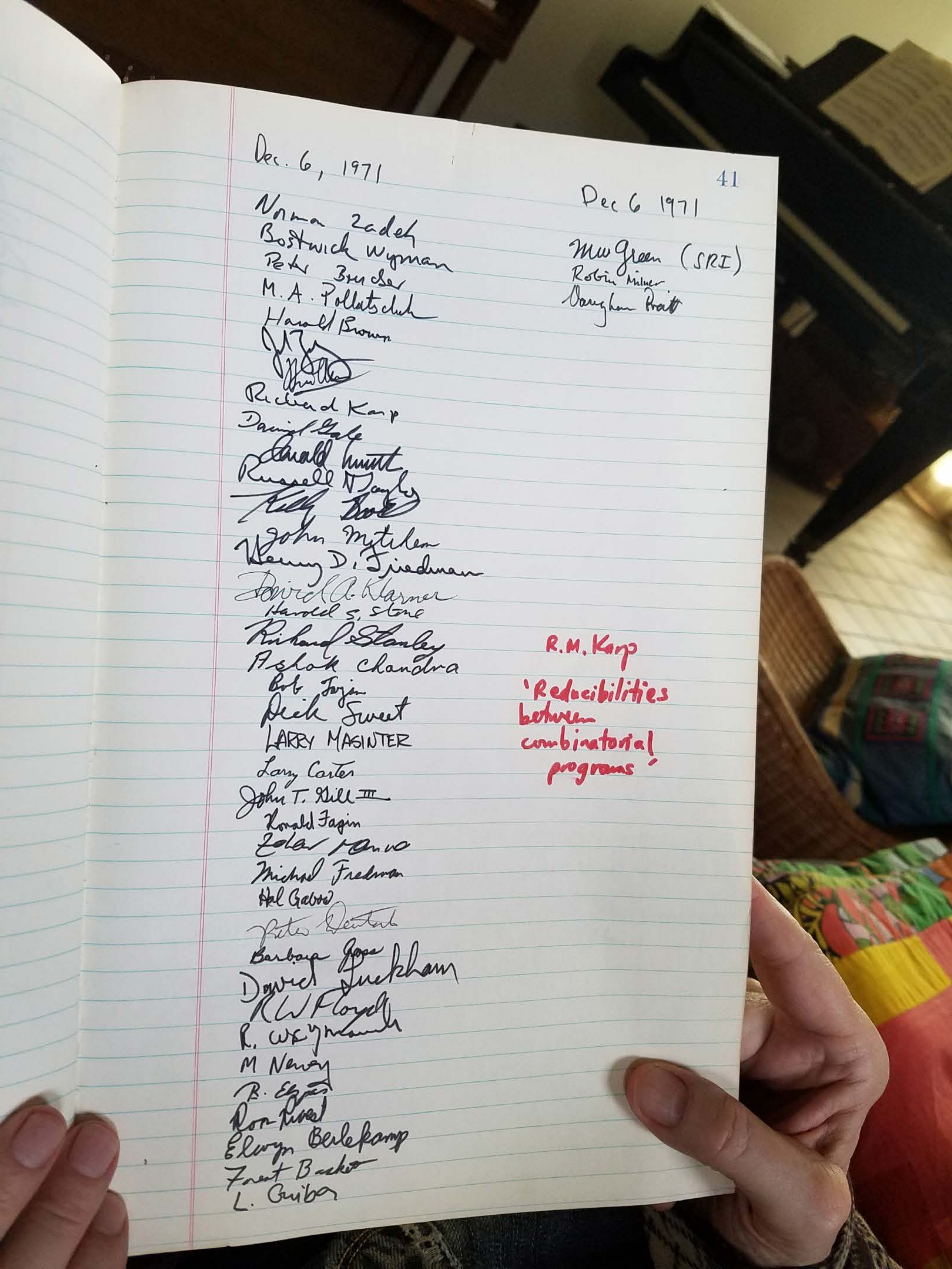}}
\caption{Knuth seminar participants}
\label{fig:karp}
\end{figure}

   After graduating from Harvard in 1970 I was an Instructor at
   \mt\ for one year before becoming a Miller Research Fellow at
   Berkeley for two years. Although Berkeley did not have the EAC
   ambience of \mt, there were nevertheless many interesting persons
   with whom I could interact, including Elwyn Berlekamp, David Gale,
   Derrick and Emma Lehmer, and Raphael and Julia
   Robinson.\footnote{Emma Lehmer and Julia Robinson were not
     officially affiliated with U.C.\ Berkeley, primarily due to
     anti-nepotism rules in effect at the time.} I also played
   duplicate bridge with Edwin Spanier. A mathematical highlight of my
   stay at Berkeley was regular visits, frequently with David Gale, to
   Stanford University in order to attend a combinatorics/computer
   science seminar held by Donald Knuth at his home on the Stanford
   campus. Figure~\ref{fig:karp} shows the participants for the
   December 6, 1971, talk of Richard Karp. (The person holding the
   book is Knuth.) This was the first public talk that discussed the P
   vs.\ NP problem.

% figure added 5/19/21

   In 1973 I returned to the exciting EAC atmosphere at M.I.T. Curtis
   Greene was an Assistant Professor 1971--1976, as was Joel Spencer
   1972--1975. The combinatorics Instructors who were there during
   some subset of the period 1973--1979 were Stephen Fisk, Karanbir
   Sarkaria, Kenneth Baclawski, Thomas Zaslavsky, Joni Shapiro (later
   Saj-nicole Joni), Dennis Stanton, Ira Gessel, and Jeff Kahn. More
   \mt\ graduate students were becoming interested in EAC. The first
   person to become my student was Emden Gansner, followed shortly
   thereafter by Ira Gessel, though Gessel was my first student to
   graduate. I had three students who graduated in the 1970's: Ira
   Gessel (1977), Emden Gansner (1978), and Bruce Sagan (1979). Paul
   Edelman, Robert Proctor, and Jim Walker were also my students
   mostly in the late 1970's, though they graduated after 1979. Walker
   was a dream student with regard to how much effort I needed to put
   in. I knew him as a graduate student for several years, but he
   talked to me about his research only occasionally and never asked
   about becoming my student. One day he walked into my office and
   asked whether some work he had written up was sufficient for a
   thesis. I looked it over for a few days and saw that it would make
   a fine thesis! Today such a scenario is not possible since the
   M.I.T.\ Math Department has instituted some strict rules for
   keeping track of graduate student progress.

   There was a combinatorics seminar run by Rota (and later me) that
   met Wednesday afternoons. Later it was expanded to Wednesdays and
   Fridays. Since \mt\ was a combinatorial magnet we had no problem
   attracting good speakers. Figure~\ref{fig:speakers} shows a list of
   some seminar speakers for the period fall 1968--fall 1970. (George
   Andrews was a Visiting Professor at the \mt\ Department of
   Mathematics for the 1970--71 academic year, thus explaining his
   many talks.) After the seminar we would frequently go to a
   student-run pub in a nearby building (Walker Memorial Hall), and
   then often to dinner. For a while (I don't recall the precise
   dates) Rota held a seminar on classical invariant theory and other
   topics entitled ``Syzygy Street.''\footnote{A \emph{syzygy} in this
     context is a certain relation among invariant polynomials. For
     readers not so familiar with American culture, the name ``Syzygy
     Street'' was inspired by the television program for preschool
     children called ``Sesame Street.''}

 \begin{figure}
\centering
\begin{tabular}{lll}
\ \ speaker & \ \ \ \ \qquad title & \ \ \quad date\\ \hline
B. Rothschild & Ramsey-type theorems & October 2, 1968\\
J. Goldman & Finite vector spaces & October 30\\
D. Kleitman & Combinatorics and statistical mechanics & unknown\\
E. Lieb & Ice is nice & January 8, 1969\\
P. O'Neil & Random 0-1 matrices & February 12\\
D. Kleitman and\\ \ \ B. Rothschild & Asymptotic enumeration of finite
topologies & March 13\\
D. Kleitman & Asymptotic enumeration of tournaments & April 22\\
E. Berlekamp & Finite Riemann Hypothesis and\\
 & \ \ error-correcting codes & April 28\\
G. Katona & Sperner-type theorems & April 29\\
E. Bender & Plane partitions and Young tableaux & April 30\\
(unknown) & Tournaments & October 6\\
S. Sherman & Monotonicity and ferromagnetism & October 20\\
G.-C.\ Rota & Exterior algebra I & December 3\\
M. Aigner & Segments of ordered sets & December 8\\
G.-C.\ Rota & Exterior algebra II & December 10\\
A. Gleason & Segments of ordered sets & December 15\\
J. Spencer & Scrambling sets & February 9, 1970\\
M. Sch\"utzenberger & Planar graphs and symmetric groups & March 4\\
D. Kleitman & Antichains in ordered sets & March 16\\
G. Szekeres & Skew block designs & April 8\\
E. Lieb & Ising and dimer problems & April 27\\
P. Erd\H{o}s & Problems of combinatorial analysis & June 3\\
D. Kleitman & (unknown) & September 29\\
G. Andrews & A partition problem of Adler & October 6\\
G. Andrews & Partitions from Euler to Gauss & October 9\\
G. Gallavotti & Some graphical enumeration problems\\
& \ \ motivated by statistical mechanics & October 13\\
G. Andrews & The Rogers-Ramanujan identities & October 16\\
R. Reid & Tutte representability and Segre arcs\\
& \ \ and caps & October 20\\
G. Andrews & Proof of Gordon's theorem & November 6\\
G. Andrews & Schur's partition theorem & November 13\\
G. Andrews & Extensions of Schur's theorem & November 20\\
M. Krieger & Some problems and conjectures & (unknown)\\
S. Fisk & Triangulations of spheres & December 1\\
D. Kleitman & Some network problems & December 8\\
\end{tabular}
\caption{\mt\ Combinatorics Seminar, 1968--1970}
\label{fig:speakers}
\end{figure}

Rota taught the course ``18.17 Combinatorial Analysis'' at \mt\ in the
fall of 1962, probably the first course on combinatorics at
\mt\footnote{It was certainly the first combinatorics course that was
  taught during or after the fall of 1958. Rota began teaching at
  \mt\ in the fall of 1959.} The listed textbooks for this course in
the \mt\ course catalog were Ore \cite{ore} and Riordan
\cite{riordan}, though much, if not all, of the material was prepared
by Rota and written up as course notes \cite{f-l-s}. These notes have
a curious feature. They were written up by G. Feldman, J. Levinger,
and Richard Stanley. However, that Richard Stanley was not me! In
fact, I was a freshman at Caltech at the time. This other Richard
Stanley (whose middle name unfortunately was John, not Peter) received
a Ph.D.\ in linguistics from \mt\ in 1969. His thesis \cite{navaho}
did have some combinatorial flavor.

\textsc{Acknowledgement.} I am grateful to Adriano Garsia, Ming-chang
Kang, Igor Pak, Victor Reiner, Christophe Reutenauer and especially
Curtis Greene for many helpful suggestions.

%%%%%%%%%%%%%%%%%%%%%%%%%%%%%%%%%%%%%%%%%%%%%%%%%%%%%%%%%%%%%%%%%

\end{document}